\documentclass[12pt]{article}
\usepackage{amssymb,amsmath}
\usepackage{cases}
\usepackage{amsfonts}
\usepackage{color,xcolor}
\usepackage[left=2.0cm,right=2.0cm,top=2.0cm,bottom=2.0cm]{geometry}
\usepackage[colorlinks,citecolor=blue,urlcolor=blue]{hyperref}

\newtheorem{theorem}{Theorem}[section]

\newtheorem{lemma}{Lemma}[section]
\newtheorem{proposition}{Proposition}[section]

\newtheorem{definition}{Definition}[section]
\newtheorem{remark}{Remark}[section]

\newcommand{\bal}{\begin{align}}
\newcommand{\bbal}{\begin{align*}}
\newcommand{\beq}{\begin{equation}}
\newcommand{\eeq}{\end{equation}}
\newcommand{\bca}{\begin{cases}}
\newcommand{\eca}{\end{cases}}

\newcommand{\pa}{\partial}
\newcommand{\fr}{\frac}

\newcommand{\De}{\Delta}

\newcommand{\dd}{\mathrm{d}}

\newcommand{\R}{\mathbb{R}}
\newcommand{\T}{\mathbb{T}}
\newcommand{\Z}{\mathbb{Z}}

\newcommand{\les}{\lesssim}

\newcommand{\bi}{\Big}
\newcommand{\g}{\big}
\begin{document}
\title{Non-uniform dependence on initial data for the Camassa-Holm equation in Besov spaces}

\author{Jinlu Li$^{1}$, Yanghai Yu$^{2}$ and Weipeng Zhu$^{3}$\footnote{E-mail: lijl29@gnnu.edu.cn; yuyanghai214@sina.com(Corresponding author); mathzwp2010@163.com}\\
\small $^1$\it School of Mathematics and Computer Sciences, Gannan Normal University, Ganzhou 341000, China\\
\small $^2$\it School of Mathematics and Statistics, Anhui Normal University, Wuhu, Anhui, 241000, China\\
\small $^3$\it School of Mathematics and Information Science, Guangzhou University, Guangzhou 510006, China}

\date{}

\maketitle\noindent{\hrulefill}

{\bf Abstract:} In the paper, we consider the initial value problem to the Camassa-Holm equation in the real-line case. Based on the local well-posedness
result and the lifespan, we proved that the data-to-solution map of this problem is not uniformly continuous in nonhomogeneous Besov spaces in the sense of Hadamard. Our obtained result improves considerably the result in \cite{H-K}.


\vskip0mm\noindent{\hrulefill}

\section{Introduction}

In what follows we are concerned with the Cauchy problem for the classical Camassa-Holm (CH) equation
\begin{equation}\label{0}
\begin{cases}
u_t-u_{xxt}+3uu_x=2u_xu_{xx}+uu_{xxx}, \; &(x,t)\in \R\times\R^+,\\
u(x,t=0)=u_0,\; &x\in \R.
\end{cases}
\end{equation}

The CH equation was firstly proposed in the context of hereditary symmetries studied in \cite{Fokas} and then was derived explicitly as a water wave equation by Camassa--Holm \cite{Camassa}. The CH equation is completely integrable \cite{Camassa,Constantin-P} with a bi-Hamiltonian structure \cite{Constantin-E,Fokas} and infinitely many conservation laws \cite{Camassa,Fokas}. Also, it admits exact peaked
soliton solutions (peakons) of the form $ce^{-|x-ct|}$ with $c>0$, which are orbitally stable \cite{Constantin.Strauss} and models wave breaking (i.e., the solution remains bounded, while its slope becomes unbounded in finite
time \cite{Constantin2000,Constantin.Escher2,Constantin.Escher3}. It is worth mentioning that the peaked solitons present the characteristic for the travelling water waves of greatest height and largest amplitude and arise as solutions to the free-boundary problem for incompressible Euler equations over a flat bed, see Refs. \cite{Constantin2,Constantin.Escher4,Constantin.Escher5,Toland} for the details. Because of the mentioned interesting and remarkable features,
 the CH equation has attracted much attention as a class of integrable shallow water wave equations in recent twenty years. Concerning the local well-posedness and ill-posedness for the Cauchy problem of the CH equation in Sobolev spaces and Besov spaces, we refer to \cite{Constantin.Escher,Constantin.Escher2,d1,d3,Guo-Yin,H-H,Li-Yin1,Ro} and the references therein. It was shown that there exist global strong solutions to the CH equation \cite{Constantin,Constantin.Escher,Constantin.Escher2} and finite time blow-up strong solutions to the CH equation \cite{Constantin,Constantin.Escher,Constantin.Escher2,Constantin.Escher3}. The existence and uniqueness of global weak solutions to the CH equation were proved in \cite{Constantin.Molinet, Xin.Z.P}. Bressan--Constantin proved the existence of the global conservative solutions \cite{Bre1} and global dissipative solutions \cite{Bre2} in $H^1(\R)$.

After the phenomenon of non-uniform continuity for some dispersive equations was studied by Kenig et.al \cite{Kenig2001}, the issue of non-uniform dependence on the initial data has been the subject of many papers. Himonas--Misio{\l}ek \cite{H-M} obtained the first result on the non-uniform dependence for the CH equation in $H^s(\T)$ with $s\geq2$ using explicitly constructed travelling wave solutions, which was sharpened to $s>\fr32$ by Himonas--Kenig \cite{H-K} on the real-line and Himonas--Kenig--Misio{\l}ek \cite{H-K-M} on the circle. It should be mentioned that Danchin \cite{d1,d3} proved the local existence and uniqueness of strong solutions to the CH equation with initial data in $B^s_{p,r}$ for $s>\max\g\{1+\fr1p, \fr32\g\}$ and $B^{\frac32}_{2,1}$. For the continuity of the solution map of the CH equation with respect to the initial data, it was proved by Li and Yin \cite{Li-Yin1}. Up to now, to our best knowledge, there is no paper concerning the non-uniform dependence on initial data for the one dimension CH euqation under the framework of Besov spaces, which is we shall investigate in this paper.

Before stating our main result, we transform the CH equation \eqref{0} equivalently into the following transport type equation
\begin{equation}\label{CH}
\begin{cases}
\partial_tu+u\pa_xu=\mathbf{P}(u), \; &(x,t)\in \R\times\R^+,\\
u(x,t=0)=u_0,\; &x\in \R,
\end{cases}
\end{equation}
where
\begin{equation}\label{CH1}
\mathbf{P}(u)=-\pa_x\g(1-\pa^2_x\g)^{-1}\bi(u^2+\fr12(\pa_xu)^2\bi).
\end{equation}
%\mathbf{P}(u)=-\pa_x(1-\pa^2_x)^{-1}\bi(\frac32u^2\bi).
Our main result is as follows.
\begin{theorem}\label{th-camassa}
Assume that $(s,p,r)$ satisfies
\begin{align}\label{con-camassa}
s>\max\bi\{1+\frac{1}{p},\frac32\bi\}, \; (p,r)\in [1,\infty]\times[1,\infty).
\end{align} Then the CH equation \eqref{CH}--\eqref{CH1} is not uniformly continuous from any bounded subset in $B^s_{p,r}$ into $\mathcal{C}([0,T];B^s_{p,r})$. More precisely, there exists two sequences of solutions $\mathbf{S}_t(f_n+g_n)$ and $\mathbf{S}_t(f_n)$ such that
\bbal
&||f_n||_{B^s_{p,r}}\lesssim 1 \quad\text{and}\quad \lim_{n\rightarrow \infty}||g_n||_{B^s_{p,r}}= 0
\end{align*}
but
\bbal
\liminf_{n\rightarrow \infty}||\mathbf{S}_t(f_n+g_n)-\mathbf{S}_t(f_n)||_{B^s_{p,r}}\gtrsim t,  \quad \forall \;t\in[0,T].
\end{align*}
\end{theorem}
\begin{remark}\label{re1}
The methods we used in proving the Theorem \ref{th-camassa} are very general and can be applied equally well to other related systems,
such as the Degasperis-Procesi equation.
\end{remark}

Another well-known integrable equation admitting peakons is the Degasperis-Procesi (DP) equation \cite{Constantin.Lannes,DP}
\begin{equation}\label{DP}
\begin{cases}
\partial_tu+u\pa_xu=-\frac32\pa_x(1-\pa^2_x)^{-1}u^2, \; &(x,t)\in \R\times\R^+,\\
u(x,t=0)=u_0,\; &x\in \R.
\end{cases}
\end{equation}
Similarly, we also have
\begin{theorem}\label{th2}
Assume that $(s,p,r)$ satisfies
\begin{align*}
s>1+\frac1p, \ p\in [1,\infty], \ r\in [1,\infty) \quad \mbox{or} \quad
s=\frac{1}{p}+1, \ p\in [1,\infty),\  r=1.
\end{align*} Then the the DP equation \eqref{DP} is not uniformly continuous from any bounded subset in $B^s_{p,r}$ into $\mathcal{C}([0,T];B^s_{p,r})$.
\end{theorem}
\begin{remark}\label{re2}
Following the procedure in the proof of Theorem \ref{th-camassa} with suitable modification, we can prove Theorem \ref{th2}. Here we will omit the details and leave it to the interested readers.
\end{remark}

\noindent\textbf{Organization of our paper.}
In Section 2, we list some notations and known results which will be used in the sequel. In Section 3, we present the local well-posedness result and establish some technical lemmas. In Section 4, we prove our main theorem. Here we give a overview of the strategy:
\begin{itemize}
  \item Choosing a sequence of approximate initial data $f_n$, which can approximate to the solution $\mathbf{S}_t(f_n)$;
  \item Considering the initial data $u^n_0=f_n+g_n$ (see Section 3.2 for the constructions of $f_n$ and $g_n$), we shall use a completely new idea. Let us make it more precise: set
  \bbal
\mathbf{w}_n=\mathbf{S}_{t}(u^n_0)-u^n_0-t\mathbf{v}_0^{n}\quad\mbox{ with }\;\mathbf{v}^n_0=-u^n_0\pa_x u^n_0,
\end{align*}
based on the special choice of $f_n$ and $g_n$, we make an important observation that the appearance of $g_n\partial_xf_n$ plays an essential role since it would not small when $n$ is large enough;
  \item The key step is to compute the error $\mathbf{w}_n$ and estimate the $B_{p,r}^s$-norm of this error;
  \item With the approximate solutions $\mathbf{S}_{t}(f_n)$ and $\mathbf{S}_{t}(u^n_0)$ were constructed, combining the precious steps, we can conclude that their distance at the initial time is converging to zero, while at any later time it is bounded below by a positive constant, namely,
\bbal
\lim_{n\rightarrow \infty}||f_n+g_n-f_n||_{B^s_{p,r}}=0
\end{align*}
but
\bbal
\liminf_{n\rightarrow \infty}||\mathbf{S}_t(f_n+g_n)-\mathbf{S}_t(f_n)||_{B^s_{p,r}}\gtrsim t\quad\text{for small }t.
\end{align*}
 That means the solution map is not uniformly continuous.
\end{itemize}

\section{Littlewood-Paley analysis}
We will use the following notations throughout this paper.
\begin{itemize}
  \item Given a Banach space $X$, we denote its norm by $\|\cdot\|_{X}$.
  \item
\noindent The symbol $A\lesssim B$ means that there is a uniform positive constant $c$ independent of $A$ and $B$ such that $A\leq cB$.
  \item Let us recall that for all $u\in \mathcal{S}'$, the Fourier transform $\mathcal{F}u$, also denoted by $\hat{u}$, is defined by
$$
\mathcal{F}u(\xi)=\hat{u}(\xi)=\int_{\R}e^{-ix\xi}u(x)\dd x \quad\text{for any}\; \xi\in\R.
$$
  \item The inverse Fourier transform allows us to recover $u$ from $\hat{u}$:
$$
u(x)=\mathcal{F}^{-1}\hat{u}(x)=\frac{1}{2\pi}\int_{\R}e^{ix\xi}\hat{u}(\xi)\dd\xi.
$$
\end{itemize}
Next, we will recall some facts about the Littlewood-Paley decomposition, the nonhomogeneous Besov spaces and their some useful properties (see \cite{B.C.D} for more details).

There exists a couple of smooth functions $(\chi,\varphi)$ valued in $[0,1]$, such that $\chi$ is supported in the ball $\mathcal{B}\triangleq \{\xi\in\mathbb{R}:|\xi|\leq \frac 4 3\}$, and $\varphi$ is supported in the ring $\mathcal{C}\triangleq \{\xi\in\mathbb{R}:\frac 3 4\leq|\xi|\leq \frac 8 3\}$. Moreover,
\begin{eqnarray*}
\chi(\xi)+\sum_{j\geq0}\varphi(2^{-j}\xi)=1 \quad \mbox{ for any } \xi\in \R.
\end{eqnarray*}
It is easy to show that $\varphi\equiv 1$ for $\frac43\leq |\xi|\leq \frac32$.

For every $u\in \mathcal{S'}(\mathbb{R})$, the inhomogeneous dyadic blocks ${\Delta}_j$ are defined as follows
\begin{numcases}{\Delta_ju=}
0, & if $j\leq-2$;\nonumber\\
\chi(D)u=\mathcal{F}^{-1}(\chi \mathcal{F}u), & if $j=-1$;\nonumber\\
\varphi(2^{-j}D)u=\mathcal{F}^{-1}\g(\varphi(2^{-j}\cdot)\mathcal{F}u\g), & if $j\geq0$.\nonumber
\end{numcases}
In the inhomogeneous case, the following Littlewood-Paley decomposition makes sense
$$
u=\sum_{j\geq-1}{\Delta}_ju\quad \text{for any}\;u\in \mathcal{S'}(\mathbb{R}).
$$

\begin{definition}(\cite{B.C.D})
Let $s\in\mathbb{R}$ and $(p,r)\in[1, \infty]^2$. The nonhomogeneous Besov space $B^s_{p,r}(\R)$ consists of all tempered distribution $u$ such that
\begin{align*}
||u||_{B^s_{p,r}(\R)}\triangleq \Big|\Big|\g(2^{js}||\Delta_j{u}||_{L^p(\R)}\g)_{j\in \Z}\Big|\Big|_{\ell^r(\Z)}<\infty.
\end{align*}
\end{definition}
\begin{remark}\label{re3}
It should be emphasized that the following embedding will be often used implicity:
$$B^s_{p,q}(\R)\hookrightarrow B^t_{p,r}(\R)\quad\text{for}\;s>t\quad\text{or}\quad s=t,1\leq q\leq r\leq\infty.$$
\end{remark}
Finally, we give some important properties which will be also often used throughout the paper.
\begin{lemma}\label{le1}(\cite{B.C.D})
Let $(p,r)\in[1, \infty]^2$ and $s>\max\bi\{1+\frac1p,\frac32\bi\}$. Then we have
\bbal
&||uv||_{B^{s-2}_{p,r}(\R)}\leq C||u||_{B^{s-2}_{p,r}(\R)}||v||_{B^{s-1}_{p,r}(\R)}.
\end{align*}
Hence, for the terms $\mathbf{P}(u)$ and $\mathbf{P}(v)$, we have
\bbal
&||\mathbf{P}(u)-\mathbf{P}(v)||_{B^{s-1}_{p,r}(\R)}\leq C||u-v||_{B^{s-1}_{p,r}(\R)}||u+v||_{B^{s}_{p,r}(\R)}.
\end{align*}
\end{lemma}

\begin{lemma}\label{le2}(\cite{B.C.D})
For $(p,r)\in[1, \infty]^2$ and $s>0$, $B^s_{p,r}(\R)\cap L^\infty(\R)$ is an algebra. Moreover, $B^{\frac{1}{p}}_{p,1}(\R)\hookrightarrow L^\infty(\R)$, and for any $u,v \in B^s_{p,r}(\R)\cap L^\infty(\R)$, we have
\bbal
&||uv||_{B^{s}_{p,r}(\R)}\leq C(||u||_{B^{s}_{p,r}(\R)}||v||_{L^\infty(\R)}+||v||_{B^{s}_{p,r}(\R)}||u||_{L^\infty(\R)}).
\end{align*}
\end{lemma}

\begin{lemma}\label{le3}(\cite{B.C.D,Li-Yin2})
Let $(p,r)\in[1, \infty]^2$ and $\sigma\geq-\min\g\{\frac1p, 1-\frac1p\g\}$. Assume that $f_0\in B^\sigma_{p,r}(\R)$, $g\in L^1([0,T]; B^\sigma_{p,r}(\R))$ and
\begin{numcases}{\pa_{x}\mathbf{u}\in}
L^1([0,T]; B^{\sigma-1}_{p,r}(\R)), & if $\sigma>1+\fr1p$ \mbox{or}\ $\sigma=1+{1\over p},\; r=1$;\nonumber\\
L^1([0,T]; B^{\sigma}_{p,r}(\R)), & if $\sigma=1+\fr1p,\; r>1$;\nonumber\\
L^1([0,T]; B^{1/p}_{p,\infty}(\R)\cap L^\infty(\R)), & if $\sigma<1+\fr1p$.\nonumber
\end{numcases}
If $f\in L^\infty([0,T]; B^\sigma_{p,r}(\R))\cap \mathcal{C}([0,T]; \mathcal{S}'(\R))$ solves the following linear transport equation:
\begin{equation*}
\quad \partial_t f+\mathbf{u}\pa_xf=g,\quad \; f|_{t=0} =f_0.
\end{equation*}

\begin{enumerate}
\item Then there exists a constant $C=C(p,r,\sigma)$ such that the following statement holds
\begin{equation*}
||f(t)||_{B^\sigma_{p,r}(\R)}\leq e^{CV(t)} \Big(||f_0||_{B^\sigma_{p,r}(\R)}+\int_0^t e^{-CV(\tau)} ||g(\tau)||_{B^\sigma_{p,r}(\R)}\mathrm{d}\tau\Big),
\end{equation*}
where
\begin{numcases}{V(t)=}
\int_0^t ||\pa_x\mathbf{u}(\tau)||_{B^{\sigma-1}_{p,r}(\R)}\mathrm{d}\tau,\; &if $\sigma>1+{1\over p}$ \mbox{or} $\sigma=1+\fr1p,\; r=1$;\nonumber\\
\int_0^t ||\pa_x\mathbf{u}(\tau)||_{B^{\sigma}_{p,r}(\R)}\mathrm{d}\tau, \; &if $\sigma=1+{1\over p},\; r>1$;\nonumber\\
\int_0^t ||\pa_x\mathbf{u}(\tau)||_{B^{1/p}_{p,\infty}(\R)\cap L^\infty(\R)}\mathrm{d}\tau,\; &if $\sigma<1+\fr1p$.\nonumber
\end{numcases}
\item If $\sigma>0$, then there exists a constant $C=C(p,r,\sigma)$ such that the following statement holds
\begin{align*}
&||f(t)||_{B^\sigma_{p,r}(\R)}\leq ||f_0||_{B^\sigma_{p,r}(\R)}+\int_0^t||g(\tau)||_{B^\sigma_{p,r}(\R)}\mathrm{d}\tau \\& \quad \quad
+C\int^t_0\Big(||f(\tau)||_{B^\sigma_{p,r}(\R)}||\pa_x\mathbf{u}(\tau)||_{L^\infty(\R)}+||\pa_x\mathbf{u}(\tau)||_{B^{\sigma-1}_{p,r}(\R)}||\pa_x f(\tau)| |_{L^\infty(\R)}\Big)\mathrm{d}\tau.
\end{align*}
\end{enumerate}
\end{lemma}

\section{Preliminaries}
Before proceeding, we recall the following local well-posedness estimates for the actual solutions.
\subsection{Local well-posedness estimates for the actual solutions}

Let us recall the local well-posedness result for the CH equation in Besov spaces.

\begin{lemma}\cite{d1,d3}\label{le4}
Assume that $(s,p,r)$ satisfies \eqref{con-camassa} and for any initial data $u_0$ which belongs to $$B_R=\g\{\psi\in B_{p,r}^s: ||\psi||_{B^{s}_{p,r}}\leq R\g\}\quad\text{for any}\;R>0.$$ Then there exists some $T=T(R,s,p,r)>0$ such that the CH equation has a unique solution $\mathbf{S}_{t}(u_0)\in \mathcal{C}([0,T];B^s_{p,r})$. Moreover, we have
\begin{align}\label{s}
||\mathbf{S}_{t}(u_0)||_{B^s_{p,r}}\leq C||u_0||_{B_{p,r}^s}.
\end{align}
\end{lemma}

\subsection{Technical Lemmas}
Firstly, we need to introduce smooth, radial cut-off functions to localize the frequency region.

 Let $\hat{\phi}\in \mathcal{C}^\infty_0(\mathbb{R})$ be an even, real-valued and non-negative funtion on $\R$ and satify
\begin{numcases}{\hat{\phi}(x)=}
1, &if $|x|\leq \frac{1}{4}$,\nonumber\\
0, &if $|x|\geq \frac{1}{2}$.\nonumber
\end{numcases}
Next, we establish the following crucial lemmas which will be used later on.
\begin{lemma}\label{ley1} For any $p\in[1,\infty]$, then there exists a positive constant $M$ such that
\begin{align}\label{m}
\liminf_{n\rightarrow \infty}\bi\|\phi^2\cos \bi(\fr{17}{12}2^nx\bi)\bi\|_{L^p}\geq M.
\end{align}
\end{lemma}
{\bf Proof}\quad Without loss of generality, we may assume that $p\in[1,\infty)$. By the Fourier iversion formula and the Fubini thereom, we see that
$$||\phi||_{L^\infty}=\sup_{x\in\R}\frac{1}{2\pi}\Big|\int_{\R}\hat{\phi}(\xi)\cos(x\xi)\dd \xi\Big|\leq \frac{1}{2\pi}\int_{\R}\hat{\phi}(\xi)\dd \xi$$
and
$$\phi(0)=\frac{1}{2\pi}\int_{\R}\hat{\phi}(\xi)\dd \xi>0.$$
Since $\phi$ is a real-valued and continuous function on $\R$, then there exists some $\delta>0$ such that
$$\phi(x)\geq \frac{||\phi||_{L^\infty}}{2}\quad\text{ for any }  x\in B_{\delta}(0).$$
Thus, we have
\bbal
\bi\|\phi^2\cos \bi(\fr{17}{12}2^nx\bi)\bi\|^p_{L^p}&\geq \frac{\phi^2(0)}4\int^\delta_{0}\bi|\cos\bi(\fr{17}{12}2^nx\bi)\bi|^p\dd x\\
&=\frac{\delta}{4}\phi^2(0)\frac{1}{2^n\widetilde{\delta}}\int^{2^n\widetilde{\delta}}_{0}|\cos x|^p\dd x\quad\text{with}\;\widetilde{\delta}=\fr{17}{12}\delta.
\end{align*}
Combining the following simple fact
\bbal
\lim_{n\rightarrow \infty}\frac{1}{2^n\widetilde{\delta}}\int_0^{2^n\widetilde{\delta}}|\cos x|^p\dd x=\frac{1}{\pi}\int^\pi_0|\cos x|^p\dd x,
\end{align*}
thus, we obtain the desired result \eqref{m}.
\begin{lemma}\label{ley2} Let $s\in\R$. Define the high frequency function $f_n$ by
$$f_n=2^{-ns}\phi(x)\sin \bi(\frac{17}{12}2^nx\bi),\quad n\gg1.$$
Then for any $\sigma\in\R$, we have
\bal\label{y1}
||f_n||_{B^\sigma_{p,r}}\leq 2^{n(\sigma-s)}||\phi||_{L^p}.
\end{align}
\end{lemma}
{\bf Proof}\quad Easy computations give that
\bbal
\hat{f}_n=2^{-ns-1}i\bi[\hat{\phi}\bi(\xi+\frac{17}{12}2^n\bi)-\hat{\phi}\bi(\xi-\frac{17}{12}2^n\bi)\bi],
\end{align*}
which implies
\bbal
\mathrm{supp} \ \hat{f}_n\subset \Big\{\xi\in\R: \ \frac{17}{12}2^n-\fr12\leq |\xi|\leq \frac{17}{12}2^n+\fr12\Big\},
\end{align*}
then, we deduce
\begin{numcases}{\Delta_j(f_n)=}
f_n, &if $j=n$,\nonumber\\
0, &if $j\neq n$.\nonumber
\end{numcases}
Thus, the definition of the Besov space tells us that the desired result \eqref{y1}.
\begin{lemma}\label{ley3} Let $s\in\R$. Define the low frequency function $g_n$ by
$$g_n=\frac{12}{17}2^{-n}\phi(x),\quad n\gg1.$$
Then there exists a positive constant $\widetilde{M}$ such that
\bbal
\liminf_{n\rightarrow \infty}||g_n\pa_xf_n||_{B^s_{p,\infty}}\geq \widetilde{M}.
\end{align*}
\end{lemma}
{\bf Proof}\quad Notice that
\bbal
\mathrm{supp} \ \hat{g}_n\subset \Big\{\xi\in\R: \ 0\leq |\xi|\leq \fr12\Big\}.
\end{align*}
Then, we have
\bbal
\mathrm{supp}\ \widehat{g_n\pa_xf_n}\subset \Big\{\xi\in\R: \ \frac{17}{12}2^n-1\leq |\xi|\leq \frac{17}{12}2^n+1\Big\},
\end{align*}
which implies
\begin{numcases}{\Delta_j\g(g_n\pa_xf_n\g)=}
g_n\pa_xf_n, &if $j=n$,\nonumber\\
0, &if $j\neq n$.\nonumber
\end{numcases}
By the definitions of $f_n$ and $g_n$, we obtain
\bbal
||g_n\pa_xf_n||_{B^s_{p,\infty}}&=2^{ns}||\De_{n}\g(g_n\pa_xf_n\g)||_{L^p}=2^{ns}||g_n\pa_xf_n||_{L^p}
\\&=\bi\|\phi^2(x)\cos \bi(\frac{17}{12}2^nx\bi)+\frac{12}{17}2^{-n}\phi(x)\pa_x\phi(x)\sin \bi(\frac{17}{12}2^nx\bi)\bi\|_{L^p}
\\&\geq \bi\|\phi^2(x)\cos \bi(\frac{17}{12}2^nx\bi)\bi\|_{L^p}-C2^{-n}.
\end{align*}
Thus, the Lemma \ref{ley1} enables us to finish the proof of the Lemma \ref{ley3}.
\section{Non-uniform continuous dependence}
In this section, we will give the proof of Theorem \ref{th-camassa}.
Firstly, based on the special choice of $f_n$, we construct approximate solutions $\mathbf{S}_{t}(f_n)$ to CH equation, then estimate the error between approximate solutions $\mathbf{S}_{t}(f_n)$ and the initial data $f_n$.
\begin{proposition}\label{pro1}
Under the assumptions of Theorem \ref{th-camassa}, we have for $k=\pm1$
\begin{equation}\label{11}
||\mathbf{S}_{t}(f_n)||_{B^{s+k}_{p,r}}\leq C2^{kn}
\end{equation}
and
\bal\label{l2}
||\mathbf{S}_{t}(f_n)-f_n||_{B^{s}_{p,r}}\leq C2^{-\frac12n(s-\frac{3}{2})}.
\end{align}
\end{proposition}
{\bf Proof}\quad The local well-posedness result (see Lemma \ref{le4}) tells us that the approximate solution $\mathbf{S}_{t}(f_n)\in \mathcal{C}([0,T];B^s_{p,r})$ and has common lifespan $T\thickapprox1$. Moreover, there holds
\bbal
||\mathbf{S}_{t}(f_n)||_{L^\infty_T(B^s_{p,r})}\leq C.
\end{align*}

By Lemmas \ref{le1}--\ref{le3}, we have for any $t\in[0,T]$ and for $k=\pm1$
\bbal
 ||\mathbf{S}_{t}(f_n)||_{B^{s+k}_{p,r}}&\leq ||f_n||_{B^{s+k}_{p,r}}+\int^t_0||\mathbf{P}\g(\mathbf{S}_{\tau}(f_n)\g)||_{B^{s+k}_{p,r}}\dd \tau+\int^t_0||\mathbf{S}_{\tau}(f_n)||_{B^{s+k}_{p,r}}||\mathbf{S}_{\tau}(f_n)||_{B^s_{p,r}}
\dd \tau\\&\leq ||f_n||_{B^{s+k}_{p,r}}+\int^t_0||\mathbf{S}_{\tau}(f_n)||_{B^{s+k}_{p,r}}||\mathbf{S}_{\tau}(f_n)||_{B^s_{p,r}}\dd \tau,
\end{align*}
which follows from Gronwall's inequality and \eqref{s} that
\bal\label{y2}
||\mathbf{S}_{t}(f_n)||_{B^{s-1}_{p,r}}\leq C2^{-n} \quad\text{and}\quad ||\mathbf{S}_{t}(f_n)||_{B^{s+1}_{p,r}}\leq C2^{n}.
\end{align}
Setting $\mathbf{\widetilde{u}}=\mathbf{S}_{t}(f_n)-f_n$, then we deduce from \eqref{CH} that
\bbal
\pa_t\mathbf{\widetilde{u}}+\mathbf{S}_{t}(f_n)\pa_x\mathbf{\widetilde{u}}=-\mathbf{\widetilde{u}}\pa_xf_n-f_n\pa_xf_n+\g[\mathbf{P}\g(\mathbf{S}_{t}(f_n)\g)-\mathbf{P}(f_n)\g]+\mathbf{P}(f_n),\quad \mathbf{\widetilde{u}}_0=0.
\end{align*}

Utilizing the Lemma \ref{le3} yields
\bal\label{4.11}
e^{-C\mathbf{V}(t)}||\mathbf{\widetilde{u}}||_{B^{s-1}_{p,r}}\les \int^t_0e^{-C\mathbf{V}(\tau)}\g\|\mathbf{\widetilde{u}}\pa_xf_n,\mathbf{P}\g(\mathbf{S}_{t}(f_n)\g)-\mathbf{P}(f_n)\g\|_{B^{s-1}_{p,r}}\dd \tau
+t||f_n\pa_xf_n,\mathbf{P}(f_n)||_{B^{s-1}_{p,r}}
\end{align}
where we denote $\mathbf{V}(t)=\int^t_0||\mathbf{S}_{t}(f_n)||_{B^{s}_{p,r}}\dd \tau$.

Combining Lemmas \ref{le1}-\ref{le2} and Lemma \ref{ley2} yields
\bbal
&||\mathbf{\widetilde{u}}\pa_xf_n||_{B^{s-1}_{p,r}}\leq C||\mathbf{\widetilde{u}}||_{B^{s-1}_{p,r}}||f_n||_{B^s_{p,r}},\\
&||\mathbf{P}\g(\mathbf{S}_{t}(f_n)\g)-\mathbf{P}(f_n)||_{B^{s-1}_{p,r}}\leq C||\mathbf{\widetilde{u}}||_{B^{s-1}_{p,r}}||\mathbf{S}_{t}(f_n),f_n||_{B^s_{p,r}},\\
&||f_n\pa_xf_n||_{B^{s-1}_{p,r}}\leq ||f_n||_{L^\infty}||f_n||_{B^s_{p,r}}+||\pa_xf_n||_{L^\infty}||f_n||_{B^{s-1}_{p,r}}\leq C2^{-sn},\\
&||\mathbf{P}(f_n)||_{B^{s-1}_{p,r}}\leq C||\mathbf{P}(f_n)||_{B^{s-\frac12}_{p,r}}\leq C2^{n(s-\frac{3}{2})}||f_n,\pa_xf_n||_{L^\infty}||f_n,\pa_xf_n||_{L^p}\leq C2^{(\frac{1}{2}-s)n}.
\end{align*}
Plugging the above inequalities into \eqref{4.11}, then by the Gronwall inequality and \eqref{s}, we infer
\bal\label{4.12}
||\mathbf{S}_{t}(f_n)-f_n||_{B^{s-1}_{p,r}}\leq C2^{(\frac{1}{2}-s)n}.
\end{align}
Applying the interpolation inequality, we obtain from \eqref{y1}, \eqref{y2} and \eqref{4.12}
\bbal
||\mathbf{S}_{t}(f_n)-f_n||_{B^{s}_{p,r}}&\leq ||\mathbf{S}_{t}(f_n)-f_n||^\frac12_{B^{s-1}_{p,r}}
||\mathbf{S}_{t}(f_n)-f_n||^\frac12_{B^{s+1}_{p,r}}\leq C2^{-\frac12n(s-\frac32)}.
\end{align*}
Thus we have finished the proof of Proposition \ref{pro1}.

To obtain the non-uniformly continuous dependence property for the CH equation, we need to construct a sequence of initial data $u^n_0=f_n+g_n$, which can not approximate to the solution $\mathbf{S}_T(u^n_0)$.
\begin{proposition}\label{pro2}
Under the assumptions of Theorem \ref{th-camassa}, we have
\bal\label{l4}
||\mathbf{S}_{t}(u^n_0)-u^n_0-t\mathbf{v}_0^{n}||_{B^{s}_{p,r}}\leq Ct^{2}+C2^{-n\min\{s-\frac{3}{2},1\}},
\end{align}
where we denote $\mathbf{v}^n_0=-u^n_0\pa_x u^n_0.$
\end{proposition}
{\bf Proof}\quad Obviously, we obtain from Lemmas \ref{ley2}--\ref{ley3} that
\bbal
||u^n_0||_{B^{s+k}_{p,r}}\leq C2^{kn} \quad \text{for }\; k\in\{0,\pm1,2\}.
\end{align*}
Then, Proposition \ref{pro1} directly tells us that for $k=\pm1$
\bal\label{13}
||\mathbf{S}_{t}(u^n_0)||_{B^{s+k}_{p,r}}\leq C2^{kn}.
\end{align}

Next, we can rewrite the solution $\mathbf{S}_{t}(u^n_0)$ as follows:
\bbal
\mathbf{S}_{t}(u^n_0)=u^n_0+t\mathbf{v}_0^{n}+\mathbf{w}_n\quad\mbox{ with }\;\mathbf{v}^n_0=-u^n_0\pa_x u^n_0.
\end{align*}
Using Lemma \ref{le1} and the fact that $B^{s-1}_{p,r}(\R)\hookrightarrow L^\infty(\R)$, we have
\bbal
||\mathbf{v}^n_0||_{B^{s-1}_{p,r}}&\leq C||u^n_0||_{B^{s-1}_{p,r}}||u^n_0||_{B^{s}_{p,r}}\leq C2^{-n},\\
||\mathbf{v}^n_0||_{B^{s+1}_{p,r}}&\leq ||u^n_0||_{L^\infty}||u^n_0||_{B^{s+2}_{p,r}}+||\pa_xu^n_0||_{L^\infty}||u^n_0||_{B^{s+1}_{p,r}}\\
&\leq C2^{-n}2^{2n}+C2^n\leq C2^n.
\end{align*}
Note that $\mathbf{w}_n=\mathbf{S}_{t}(u^n_0)-u^n_0-t\mathbf{v}_0^{n}$, then we can deduce that $\mathbf{w}_n$ satisfy the following equation
\begin{eqnarray}\label{er}
\left\{\begin{array}{ll}
\pa_t\mathbf{w}_n+\mathbf{S}_{t}(u^n_0)\pa_x\mathbf{w}_n=-t\g(u^n_0\pa_x\mathbf{v}^{n}_0+\mathbf{v}^n_0\pa_xu^n_0-2\mathcal{A}_1\g)
-t^2\g(\mathbf{v}^n_0\pa_x\mathbf{v}^n_0-\mathbf{P}(\mathbf{v}^n_0)\g)\\~~~~~~~~~~~~~~~~~~~~~~~~~~~~~-\mathbf{w}_n\pa_x(u^n_0+t\mathbf{v}^n_0)+\mathcal{A}_2
+\mathcal{A}_3+\mathbf{P}(u^n_0),\\
\mathbf{w}_n(x,t=0)=0,\end{array}\right.
\end{eqnarray}
where
\bbal\mathcal{A}_1&=-\pa_x(1-\pa^2_x)^{-1}\bi(u^n_0\mathbf{v}^n_0+\fr12\partial_xu^n_0\pa_x\mathbf{v}^n_0\bi)\\
\mathcal{A}_2&=-\pa_x(1-\pa^2_x)^{-1}\bi(\mathbf{w}_n\mathbf{S}_{t}(u_n^0)+\fr12\partial_x\mathbf{w}_n\pa_x\mathbf{S}_{t}(u_n^0)\bi)\\
\mathcal{A}_3&=-\pa_x(1-\pa^2_x)^{-1}\bi(\mathbf{w}_n(u^0_n+t\mathbf{v}^0_n)+\fr12\pa_x\mathbf{w}_n\partial_x(u^0_n+t\mathbf{v}^0_n)\bi).
\end{align*}
The local well-posedness result (see Lemma \ref{ley1}) tells us that the approximate solution $\mathbf{S}_{t}(u_0^n)\in \mathcal{C}([0,T];B^s_{p,r})$ and has common lifespan $T\thickapprox1$.

Utilizing Lemma \ref{le1} to \eqref{er}, for $k\in\{-1,0\}$, we have for all $t\in[0,T]$,
\bal\label{yy}
||\mathbf{w}_n||_{B^{s+k}_{p,r}}\leq&~ C\int^t_0||\mathbf{w}_n||_{B^{s+k}_{p,r}}||u^n_{0},\mathbf{v}^n_0,\mathbf{S}_{\tau}(u^n_0)||_{B^s_{p,r}}\dd \tau+Ct||\mathbf{P}(u^n_0)||_{B^{s+k}_{p,r}}
\nonumber\\&+C(k+1)\int^t_0||\mathbf{w}_n||_{B^{s-1}_{p,r}}||u^n_{0},\mathbf{v}^n_0||_{B^{s+1}_{p,r}}\dd \tau
\nonumber\\&+Ct^2||u^n_0\pa_x\mathbf{v}_{0}^n,\mathbf{v}^n_0\pa_xu_{0}^n,\mathcal{A}_1||_{B^{s+k}_{p,r}}+Ct^3||\mathbf{v}^n_0\pa_x\mathbf{v}_{0}^n,
\mathbf{P}\big(\mathbf{v}^n_0\big)||_{B^{s+k}_{p,r}}.
\end{align}
Next, we need to estimate the above terms one by one.

{\bf \underline{Case $k=-1$}.}\quad From Lemma \ref{le1}, we have
\bbal
&||u^n_0\pa_x\mathbf{v}^n_0||_{B^{s-1}_{p,r}}\leq
C||u^n_0||_{B^{s-1}_{p,r}}||\mathbf{v}^n_0||_{B^{s}_{p,r}}\leq C2^{-n},
\\&||\mathbf{v}^n_0\pa_xu^n_0||_{B^{s-1}_{p,r}}\leq
C||\mathbf{v}^n_0||_{B^{s-1}_{p,r}}||u^n_0||_{B^{s}_{p,r}}\leq C2^{-n},
\\&||\mathbf{v}^n_0\pa_x\mathbf{v}^n_0||_{B^{s-1}_{p,r}}\leq
C||\mathbf{v}^n_0||_{B^{s-1}_{p,r}}||\mathbf{v}^n_0||_{B^{s}_{p,r}}\leq C2^{-n},
\\&||\mathcal{A}_1||_{B^{s-1}_{p,r}}\leq C||u^n_0||_{B^{s-1}_{p,r}}||\mathbf{v}^n_0||_{B^{s}_{p,r}}\leq C2^{-n},
\\&||\mathbf{P}\big(\mathbf{v}^n_0\big)||_{B^{s-1}_{p,r}}\leq C||\mathbf{v}^n_0||_{B^{s-1}_{p,r}}||\mathbf{v}^n_0||_{B^{s}_{p,r}}\leq C2^{-n}.
\end{align*}
Due to $u_0^n=f_n+g_n$, one has
$$\mathbf{P}(u^n_0)=\mathbf{P}(f_n)+\mathbf{P}(g_n)+\widetilde{\mathbf{P}}(f_n,g_n)$$
where $\widetilde{\mathbf{P}}(f_n,g_n)=-\pa_x(1-\pa^2_x)^{-1}\g(2f_ng_n+\partial_xf_n\pa_xg_n\g).$

By Lemma \ref{le1}, we have
\bbal
&||\mathbf{P}(f_n)||_{B^{s-1}_{p,r}}\leq C2^{n(s-\frac{3}{2})}||f_n,\pa_xf_n||_{L^\infty}||f_n,\pa_xf_n||_{L^p}\leq 2^{(\frac{1}{2}-s)n},
\\&||\mathbf{P}(g_n)||_{B^{s-1}_{p,r}}\leq C||g_n||^2_{B^s_{p,r}}\leq C2^{-2n},
\\&||\widetilde{\mathbf{P}}(f_n,g_n)||_{B^{s-1}_{p,r}}\leq C2^{n(s-2)}||f_n,\pa_xf_n||_{L^p}||g_n,\pa_xg_n||_{L^\infty}\leq C2^{-2n},
\end{align*}
which tell us that
\bal\label{lj}
&||\mathbf{P}(u_0^n)||_{B^{s-1}_{p,r}}\leq C2^{-n\min\{s-\frac{1}{2},2\}}.
\end{align}
Gathering all the above estimates together with \eqref{yy} and using the Gronwall inequality yields
\bbal
||\mathbf{w}_n||_{B^{s-1}_{p,r}}\leq Ct^22^{-n}+C2^{-n\min\{s-\frac{1}{2},2\}}.
\end{align*}

{\bf \underline{Case $k=0$}.}\quad From Lemmas \ref{le1}-\ref{le2}, we have
\bbal
&||u^n_0\pa_x\mathbf{v}^n_0||_{B^{s}_{p,r}}\les||u^n_0||_{B^{s-1}_{p,r}}||\mathbf{v}^n_0||_{B^{s+1}_{p,r}}+||u^n_0||_{B^{s}_{p,r}}||\mathbf{v}^n_0||_{B^{s}_{p,r}}
\les1,
\\&||\mathbf{v}^n_0\pa_xu^n_0||_{B^{s}_{p,r}}\les||\mathbf{v}^n_0||_{B^{s-1}_{p,r}}||u^n_0||_{B^{s+1}_{p,r}}+||\mathbf{v}^n_0||_{B^{s}_{p,r}}||u^n_0||_{B^{s}_{p,r}}
\les1,
\\&||\mathbf{v}^n_0\pa_x\mathbf{v}^n_0||_{B^{s}_{p,r}}\les||\mathbf{v}^n_0||_{B^{s-1}_{p,r}}||\mathbf{v}^n_0||_{B^{s+1}_{p,r}}+||\mathbf{v}^n_0||_{B^{s}_{p,r}}||\mathbf{v}^n_0||_{B^{s}_{p,r}}
\les1,
\\&||\mathcal{A}_1||_{B^{s}_{p,r}}\les||u^n_0||_{B^{s}_{p,r}}||\mathbf{v}^n_0||_{B^{s}_{p,r}}\les
1,
\\&||\mathbf{P}\big(\mathbf{v}^n_0\big)||_{B^{s}_{p,r}}\les||\mathbf{v}^n_0||_{B^{s}_{p,r}}||\mathbf{v}^n_0||_{B^{s}_{p,r}}\les1\\
&||\mathbf{P}(f_n)||_{B^{s}_{p,r}}\les2^{n(s-1)}||f_n,\pa_xf_n||_{L^\infty}||f_n,\pa_xf_n||_{L^p}\les 2^{(1-s)n},
\\&||\mathbf{P}(g_n)||_{B^{s}_{p,r}}\les||g_n||^2_{B^s_{p,r}}\leq C2^{-2n},
\\&||\widetilde{\mathbf{P}}(f_n,g_n)||_{B^{s}_{p,r}}\les2^{n(s-1)}||f_n,\pa_xf_n||_{L^p}||g_n,\pa_xg_n||_{L^\infty}\leq C2^{-n},
\end{align*}
Gathering all the above estimates together with \eqref{yy} and using the Gronwall inequality yields
\bbal
||\mathbf{w}_n||_{B^{s}_{p,r}}&\leq Ct^2+Ct2^{-n\min\{s-1,1\}}+C\int^t_02^n||\mathbf{w}_n||_{B^{s-1}_{p,r}}\dd \tau
\\&\leq Ct^2+C2^{-n\min\{s-\frac{3}{2},1\}}.
\end{align*}
Thus, we completed the proof of Proposition \ref{pro2}.

With the propositions \ref{pro1}--\ref{pro2} in hand, we can prove Theorem \ref{th-camassa}.

{\bf Proof of Theorem \ref{th-camassa}}\quad
Obviously, we have
\bbal
||u^n_0-f_n||_{B^s_{p,r}}=||g_n||_{B^s_{p,r}}\leq C2^{-n},
\end{align*}
which means that
\bbal
\lim_{n\to\infty}||u^n_0-f_n||_{B^s_{p,r}}=0.
\end{align*}
Furthermore, we deduce that
\bal\label{yyh}
\quad \ ||\mathbf{S}_{t}(u^n_0)-\mathbf{S}_{t}(f_n)||_{B^s_{p,r}}=&~||t\mathbf{v}^n_{0}+g_n+f_n-\mathbf{S}_{t}(f_n)+\mathbf{w}_n||_{B^s_{p,r}}\nonumber\\
\geq&~ ||t\mathbf{v}^n_{0}||_{B^s_{p,r}}-||g_n||_{B^s_{p,r}}-||f_n-\mathbf{S}_{t}(f_n)||_{B^s_{p,r}}-||\mathbf{w}_n||_{B^s_{p,r}}\nonumber\\
\geq&~ t||\mathbf{v}^n_0||_{B^s_{p,\infty}}-C2^{-\frac{1}{2}n\min\{s-\frac{3}{2},1\}}-Ct^{2}.
\end{align}
Notice that
$$
-\mathbf{v}^n_0=f_n\pa_xf_n+f_n\pa_xg_n+g_n\pa_xg_n+g_n\pa_xf_n,
$$
by simple calculation, we obtain
\bbal
||f_n\pa_xf_n||_{B^s_{p,r}}&\leq ||f_n||_{L^\infty}||f_n||_{B^{s+1}_{p,r}}+||\pa_xf_n||_{L^\infty}||f_n||_{B^{s}_{p,r}}\leq C2^{-n(s-1)},\\
||f_n\pa_xg_n||_{B^s_{p,r}}&\leq ||f_n||_{B^s_{p,r}}||g_n||_{B^{s+1}_{p,r}}\leq C2^{-n},\\
||g_n\pa_xg_n||_{B^s_{p,r}}&\leq ||g_n||_{B^s_{p,r}}||g_n||_{B^{s+1}_{p,r}}\leq C2^{-2n}.
\end{align*}
Hence, it follows from \eqref{yyh} and Lemma \ref{ley3} that
\bbal
\liminf_{n\rightarrow \infty}||\mathbf{S}_t(f_n+g_n)-\mathbf{S}_t(f_n)||_{B^s_{p,r}}\gtrsim t\quad\text{for} \ t \ \text{small enough}.
\end{align*}
This completes the proof of Theorem \ref{th-camassa}.

\vspace*{1em}
\noindent\textbf{Acknowledgements.}  J. Li is supported by the National Natural Science Foundation of China (Grant No.11801090). Y. Yu is supported by the Natural Science Foundation of Anhui Province (No.1908085QA05). W. Zhu is partially supported by the National Natural Science Foundation of China (Grant No.11901092) and Natural Science Foundation of Guangdong Province (No.2017A030310634).
%\vspace*{1em}


\begin{thebibliography}{99}
\linespread{0}\addtolength{\itemsep}{-1.0ex}

\bibitem{B.C.D} H. Bahouri, J. Y. Chemin and R. Danchin, \textit{Fourier Analysis and Nonlinear Partial Differential Equations}, {Grundlehren der Mathematischen Wissenschaften}, vol. 343, Springer-Verlag, Berlin, Heidelberg, 2011.
\bibitem{Bre1}A. Bressan and A. Constantin, \textit{Global conservative solutions of the Camassa-Holm equation}, {Arch. Ration. Mech. Anal.}, {\bf183} (2007), 215-239.
\bibitem{Bre2} A. Bressan and A. Constantin, \textit{Global dissipative solutions of the Camassa-Holm equation}, {Anal. Appl.}, {\bf5} (2007), 1--27.
\bibitem{Camassa}R. Camassa and D. D. Holm, \textit{An integrable shallow water equation with peaked solitons}, {Phys. Rev. Lett.}, {\bf71} (1993), 1661--1664.

\bibitem{Constantin2000} A. Constantin, \textit{Existence of permanent and breaking waves for a shallow water equation: ageometric approach}, {Ann. Inst. Fourier} {\bf 50} (2000), 321--362.

\bibitem{Constantin-E} A. Constantin, \textit{The Hamiltonian structure of the Camassa-Holm equation}, {Exposition. Math.}, {\bf15} (1997), 53--85.

\bibitem{Constantin-P} A. Constantin, \textit{On the scattering problem for the Camassa-Holm equation}, {R. Soc. Lond. Proc. Ser. A Math. Phys. Eng. Sci.}, {\bf 457} (2001), 953--970.

\bibitem{Constantin} A. Constantin, \textit{Existence of permanent and breaking waves for a shallow water equation: a geometric approach}, {Ann. Inst. Fourier (Grenoble)}, {\bf50} (2000), 321--362.

\bibitem{Constantin2} A. Constantin, \textit{The trajectories of particles in Stokes waves}, {Invent. Math.}, {\bf166} (2006), 523--535.

\bibitem{Constantin.Escher} A. Constantin and J. Escher, \textit{Global existence and blow-up for a shallow water equation}, {Ann. Scuola Norm. Sup. Pisa Cl. Sci. (4)},  {\bf26} (1998), 303--328.

\bibitem{Constantin.Escher2} A. Constantin and J. Escher, \textit{Well-posedness, global existence, and blowup phenomena for a periodic quasi-linear hyperbolic equation}, {Comm. Pure Appl. Math.}, {\bf51} (1998), 475--504.

\bibitem{Constantin.Escher3} A. Constantin and J. Escher, \textit{Wave breaking for nonlinear nonlocal shallow water equations}, {Acta Math.}, {\bf181} (1998), 229-243.

\bibitem{Constantin.Escher4} A. Constantin and J. Escher, \textit{Particle trajectories in solitary water waves}, {Bull. Amer. Math. Soc.}, {\bf44} (2007), 423--431.

\bibitem{Constantin.Escher5} A. Constantin and J. Escher, \textit{Analyticity of periodic traveling free surface water waves with vorticity}, {Ann. of Math.}, {\bf173} (2011), 559--568.

\bibitem{Constantin.Lannes} A. Constantin and D. Lannes, \textit{The hydrodynamical relevance of the Camassa-Holm and Degasperis-Procesi equations}, {Arch. Ration. Mech. Anal.}, {\bf192} (2009), 165-186.

\bibitem{Constantin.Molinet} A. Constantin and L. Molinet \textit{Global weak solutions for a shallow water equation}, {Comm. Math. Phys.}, {\bf211} (2000), 45--61.

\bibitem{Constantin.Strauss}  A. Constantin and W. A. Strauss, \textit{Stability of peakons}, {Comm. Pure Appl. Math.}, {\bf53} (2000), 603--610.

\bibitem{d1}R. Danchin, \textit{A few remarks on the Camassa-Holm equation}. {Differential Integral Equations}, {\bf14} (2001), 953--988.


\bibitem{d3} R. Danchin, \textit{A note on well-posedness for Camassa-Holm equation}, {J. Differential Equations}, {\bf192} (2003), 429--444.

\bibitem{DP} A. Degasperis and M. Procesi, \textit{Asymptotic integrability, in: Symmetry and Perturbation Theory}, Rome, 1998, World Sci. Publ., River Edge, NJ, (1999), p.23.


\bibitem{Fokas} B. Fuchssteiner and A. Fokas, \textit{Symplectic structures, their B\"{a}cklund transformation and hereditary symmetries}, {Phys. D}, {\bf 4} (1981/82), 47--66.

\bibitem{Guo-Yin} Z. Guo, X. Liu, M. Luc and Z. Yin, \textit{Ill-posedness of the Camassa-Holm and related equations in the critical space}, {J. Differential Equations}, {\bf266} (2019), 1698--1707.

\bibitem{H-K} A. Himonas and C. Kenig, \textit{Non-uniform dependence on initial data for the CH equation on the line}, {Diff. Integral Eqns}, {\bf22} (2009), 201--224.

\bibitem{H-K-M} A. Himonas, C. Kenig and Misio{\l}ek \textit{Non-uniform dependence for the periodic CH equation}, {Commun. Partial Diff. Eqns}, {\bf35} (2010), 1145--1162.

\bibitem{H-H} A. Himonas and C. Holliman, \textit{The Cauchy problem for the Novikov equation}, {Nonlinearity}, {25} (2012), 449--479.

\bibitem{H-M} A. Himonas and G. Misio{\l}ek, \textit{Non-uniform dependence on initial data of solutions to the Euler equations of hydrodynamics},  {Comm. Math. Phys}., 296 (2010), 285--301.

\bibitem{Kenig2001} C. Kenig, G. Ponce, L. Vega, \textit{On the ill-posedness of some canonical dispersive equations}, {Duke Math}., {\bf106} (2001) 617--633.

\bibitem{Li-Yin1} J. Li and Z. Yin, \textit{Remarks on the well-posedness of Camassa-Holm type equations in Besov spaces}, {J. Differential Equations}, {\bf261} (2016), 6125-6143.

\bibitem{Li-Yin2} J. Li  and Z. Yin, \textit{Well-posedness and analytic solutions of the two-component Euler-Poincar\'{e} system}, {Monatsh. Math.}, 183 (2017), 509--537.


\bibitem{Ro} G. Rodr\'{i}guez-Blanco, \textit{On the Cauchy problem for the Camassa-Holm equation}, {Nonlinear Anal.}, {\bf46} (2001), 309--327.

\bibitem{Toland} J. F. Toland, \textit{Stokes waves}, {Topol. Methods Nonlinear Anal.}, {\bf7} (1996), 1--48.



\bibitem{Xin.Z.P} Z. Xin and P. Zhang, \textit{On the weak solutions to a shallow water equation}, {Comm. Pure Appl. Math.}, {\bf53} (2000), 1411--1433.



\end{thebibliography}
\end{document}